\documentstyle[12pt]{article}
\begin{document}
\begin{center}
{\Large\bf ON FRACTIONAL KINETIC EQUATIONS}\\[1cm]
{\large R.K. SAXENA}\\
Department of Mathematics and Statistics, 
Jai Narain Vyas University Jodhpur
342001, INDIA\\[0.5cm]
{\large A.M. MATHAI}\\ 
Department of Mathematics and Statistics, McGill University,\\
805 Sherbooke Street West, Montreal, CANADA H3A 2K6\\[0.5cm]
{\large H.J. HAUBOLD}\\ 
Office for Outer Space Affairs, United Nations,\\ 
P.O. Box 500, A-1400 Vienna, AUSTRIA\\[0.5cm]
\end{center}
\noindent

{\bf Abstract.} The subject of this paper is to derive the solution of generalized fractional kinetic equations. The results are obtained in a compact form containing the Mittag-Leffler function, which naturally occurs whenever one is dealing with fractional integral equations. The results derived in this paper provide an extension of a result given by Haubold and Mathai in a recent paper (Haubold and Mathai, 2000).
\section{Introduction and Preliminaries}
In terms of Pochammer's symbol
$$(\alpha)_n=\left\{^{1,n=0}_{\alpha(\alpha+1)\ldots(\alpha+n-1),n\in N}\right.$$
we can express the binomial series as
\begin{equation}
(1-x)^{-\alpha}=\sum^\infty_{r=0}\frac{(\alpha)_rx^r}{r!}.
\end{equation}
The Mittag-Leffler function is defined by
\begin{equation}
E_\alpha(x):=\sum^\infty_{n=0}\frac{z^n}{\Gamma(\alpha n+1)},
\end{equation}
This function was defined and studied by Mittag-Leffler (Mittag-Leffler, 1902, 1905).
We note that this function is a direct generalization of an exponential function, since
$$E_1(z):=exp(z).$$
It also includes the error functions and other related functions, for we have
\begin{equation}
E_{1/2}(\pm z^{1/2})=e^z[1+erf(\pm z^{1/2})]=e^z erfc(\mp z^{1/2}),
\end{equation}
where \begin{equation}
erf(z):=\frac{2}{\pi^{1/2}}\int^z_0e^{-u^2}du, erfc(z):= 1-erf(z), z\in C.
\end{equation}
The equation
\begin{equation}
E_{\alpha,\beta}(z):=\sum^\infty_{n=0}\frac{z^n}{\Gamma(\alpha n+\beta)}
\end{equation}
gives a generalization of the Mittag-Leffler function. This generalization was studied by Wiman (1905), Agarwal (1953), Humbert (1953) and Humbert and Agarwal (1953) and several others.
When $\beta = 1$, (5) reduces to (2).
Both the functions defined by (2) and (5) are entire functions of order $1/\alpha$ and type 1. A detailed account of these functions is available from the monograph of Erd\'{e}lyi et al. (1955).
The Laplace transform of $E_{\alpha, \beta}(z)$ follows from the integral
\begin{equation}
\int^\infty_0e^{-pt}t^{\beta-1}E_{\alpha, \beta}(\lambda at^\alpha)dt=p^{-\beta}(1-ap^{-\alpha})^{-1},
\end{equation}
where $Re(p)>|a|^{1/\alpha}, Re(\beta)>0$, which can be established by means of the Laplace integral
\begin{equation}
\int^\infty_0e^{-pt}t^{\rho-1}dt=\Gamma(\rho)/p^\rho,
\end{equation}
where $Re(p)>0, Re(\rho)>0$.
The Riemann-Liouville operator of fractional integration is defined as
\begin{equation}
_aD_t^{-\nu}f(t)=\frac{1}{\Gamma(\nu)}\int^t_a  f(u)(t-u)^{\nu-1}du, \nu>0,
\end{equation}
with $_aD_t^0 f(t)=f(t)$ (Oldham and Spanier, 1974; Miller and Ross, 1993; Srivastava and Saxena, 2001).
By integrating the standard kinetic equation
\begin{equation}
\frac{d}{dt}N_i(t)=-c_iN_i(t), (c_i>0),
\end{equation}
it is derived that (Haubold and Mathai, 2000)
\begin{equation}
N_i(t)-N_0=-c_i\;\;_0D_t^{-1}N_i(t),
\end{equation}
where $_0D_t^{-1}$ is the standard Riemann integral operator. Here we recall that, in the original paper of Haubold and Mathai (2000), the number density of species $i, N_i=N_i(t)$, is a function of time and $N_i(t=0)=N_0$ is the number density of species $i$ at time $t=0$. By dropping the index $i$ in (10), the solution of its generalized form
\begin{equation}
N(t)-N_0=-c^\nu\;\;_0D_t^{-\nu}N(t),
\end{equation}
is obtained (Haubold and Mathai, 2000) as 
\begin{equation}
N(t)=N_0\sum^\infty_{k=0}\frac{(-1)^k(ct)^{\nu k}}{\Gamma(\nu k+1)},
\end{equation}
By virtue of (2) we can rewrite (12) in terms of the Mittag-Leffler function in a compact form as
\begin{equation}
N(t)=N_0 E_\nu(-c^\nu t^\nu), \nu>0.
\end{equation}
In this paper we investigate the solutions of three generalized forms of (11). The results are obtained in a compact form in terms of the generalized Mittag-Leffler function and given in the form of three theorems.
\section{Generalized Fractional Kinetic Equations}
{\bf Theorem 1.} If $\nu>0, \mu>0,$ then the solution of the integral equation
\begin{equation}
N(t)-N_0t^{\mu-1}= -c^\nu\;\; _0D_t^{-\nu}N(t),
\end{equation}
is given by
\begin{equation}
N(t)=N_0\Gamma(\mu)t^{\mu-1} E_{\nu, \mu}(-c^\nu t^\nu),
\end{equation}
where $E_{\nu, \mu}(.)$ is the generalized Mittag-Leffler function defined by (5).\\
{\it Proof.} We know that (Erd\'{e}lyi et al., 1954) the Laplace transform of the Riemann-Liouville fractional integral is given by
\begin{equation}
L\left\{_0D_t^{-\sigma} f(t); p \right\} = p^{-\sigma}F(p),
\end{equation}
where
\begin{equation}
F(p)=\int^\infty_{u=0}e^{-pu} f(u) du.
\end{equation}
Projecting the equation (14) to Laplace transform, it gives
\begin{equation}
N(p)=L\left\{N(t); p\right\} = \frac{N_0\Gamma(\mu)}{p^\mu\left\{1+(p/c)^{-\nu}\right\}}.
\end{equation}
By virtue of the relation
\begin{equation}
L^{-1}\left\{p^{-\rho}\right\}= \frac{t^{\rho-1}}{\Gamma(\rho)}, Re(\rho)>0,
\end{equation}
it is found that
\begin{eqnarray}
L^{-1}\left[\frac{N_0\Gamma(\mu)}{p^\mu\left\{1+(p/c)^{-\nu}\right\}}\right]& = & N_0\Gamma(\mu)L^{-1}\left\{p^{-\mu}\sum^\infty_{r=0}\frac{(1)_r[-(p/c)^{-\nu}]^r}{(r)!}\right\}\nonumber\\
& = & N_0\Gamma(\mu)\sum^\infty_{r=0}(-1)^rc^{\nu r}L^{-1}\left\{p^{-\mu-\nu r}\right\}\nonumber\\
& = & N_0\Gamma(\mu)\sum^\infty_{r=0}(-1)^r c^{r\nu}\frac{t^{\mu+r \nu -1}}{\Gamma(r\nu+\mu)}\nonumber\\
& = & N_0\Gamma(\mu)t^{\mu-1} E_{\nu, \mu}(-c^\nu t^\nu).
\end{eqnarray}
The result (15) now readily follows by taking inverse Laplace transform of (18). For $\mu=1$, we obtain the result given by Haubold and Mathai (2000).\\
{\bf Theorem 2.} If $\nu>0, c>0, d>0, \mu>0, Re(p)>|d|^{\nu/\alpha}, \;\;c\neq d$ then for the solution of the equation
\begin{equation}
N(t)-N_0 t^{\mu-1}E_{\nu,\mu}(-d^\nu t^\nu)=-c^\nu\;\;_0D_t^{-\nu}N(t),
\end{equation}
there holds the formula
\begin{equation}
N(t)=N_0\frac{t^{\mu-\nu-1}}{c^\nu-d^\nu}\left[E_{\nu,\mu-\nu}(-d^\nu t^\nu)-E_{\nu,\mu-\nu}(-c^\nu t^\nu)\right].
\end{equation}
{\it Proof.} Projecting (21) to Laplace transform and using (5) and (16), we obtain
\begin{eqnarray*}
N(p)=L\left\{N(t);p\right\} & = & N_0\frac{p^{-\mu}[1+(p/d)^{-\nu}]^{-1}}{[1+(p/c)^{-\nu}]} \\
& = &\frac{N_0 p^{\nu-\mu}}{c^\nu-d^\nu}\left[\sum^\infty_{r=0}\frac{(-1)^r(1)_r(d/p)^{r\nu}}{(r)!}-\sum^\infty_{r=0}\frac{(-1)^r(1)_r(c/p)^{r\nu}}{(r)!}\right]
\end{eqnarray*}
Hence
\begin{eqnarray*}
L^{-1}\left\{N(p)\right\} & = & N(t)\nonumber\\
&=&\frac{N_0}{c^\nu-d^\nu}\left[\sum^\infty_{r=0}(-1)^rd^{r\nu}L^{-1}(p^{-(\mu+r\nu-\nu)})-\sum^\infty_{r=0}(-1)^rc^{r\nu}L^{-1}(p^{-(\mu+r\nu-\nu)})\right]\nonumber\\
&=& \frac{N_0}{c^\nu-d^\nu}\left[\sum^\infty_{r=0}\frac{(-1)^rd^{r\nu}t^{\mu+r\nu-\nu-1}}{\Gamma(r\nu+\mu-\nu)}-\sum^\infty_{r=0}\frac{(-1)^rc^{r\nu}t^{\mu+r\nu-\nu-1}}{
\Gamma(r\nu+\mu-\nu)}\right]\nonumber\\
&=&\frac{N_0t^{\mu-\nu-1}}{c^\nu-d^\nu}\left[E_{\nu,\mu-\nu}(-d^\nu t^\nu)-E_{\nu,\mu-\nu}(-c^\nu t^\nu)\right]
\end{eqnarray*}
This completes the proof of (22).\\
When $\mu=\nu+1$, theorem 2 reduces to\\
{\bf Corollary 2.1.} If $\nu>0,c>0, d>0, c\neq d, Re(p)>|d|^{\nu/\alpha}$, then for the solution of
\begin{equation}
N(t)-N_0t^\nu E_{\nu, \nu+1}(-d^\nu t^\nu)= - c^\nu\;\; _0D_t^{-\nu} N(t),
\end{equation}
the  following result holds
\begin{equation}
N(t)=\frac{N_0}{c^\nu-d^\nu}\left[E_\nu(-d^\nu t^\nu)-E_\nu(-c^\nu t^\nu)\right].
\end{equation}
On the other hand if $d \rightarrow 0$ in (22), we arrive at\\
{\bf Corollary 2.2.} If $c>0, \nu>0, \mu>0, Re(p)>|d|^{\nu/\alpha}$, then for the solution of
\begin{equation}
N(t)- \frac{N_0 t^{\mu-1}}{\Gamma(\mu)}=-c^\nu\;\;_0D_t^{-\nu} N(t),
\end{equation}
the following result holds
\begin{equation}
N(t)=\frac{N_0 t^{\mu-\nu-1}}{c^\nu}\left[\frac{1}{\Gamma(\mu-\nu)}-E_{\nu, \mu-\nu}(-c^\nu t^\nu)\right].\\
\end{equation}
{\it Note:} When $\mu=\nu+1$, then for the solution of
\begin{equation}
N(t)-N_0\frac{t^\nu}{\Gamma(\nu+1)}=-c^\nu\;\; _0D_t^{-\nu}N(t),
\end{equation}
there holds the formula
\begin{equation}
N(t)=\frac{N_0}{c^\nu}\left[1-E_\nu(-c^\nu t^\nu)\right],
\end{equation}
where $c>0, \nu>0$.\\
The case $c=d$ is given by\\
{\bf Theorem 3.} If $c>0, \nu>0, \mu>0$, then for the solution of the equation
\begin{equation}
N(t)-N_0t^{\mu-1} E_{\nu, \mu}(-c^\nu t^\nu)= -c^\nu\;\; _0D_t^{-\nu},
\end{equation}
the following result holds
\begin{equation}
N(t)=\frac{N_0}{\nu}t^{\mu-1}\left[E_{\nu, \mu-1}(-c^\nu t^\nu)+(1+\nu-\mu)E_{\nu, \mu}(-c^\nu t^\nu)\right].
\end{equation}
{\it Proof.} Proceeding in a similar manner, it is observed that
\begin{eqnarray*}
N(t)& = & N_0L^{-1}\left\{p^{-\mu}(1+c^\nu p^{-\nu})^{-2}\right\} = N_0\sum^\infty_{r=0}\frac{(-1)^r(2)_rc^{r\nu}L^{-1}\left\{p^{-\mu-r\nu}\right\}}{(r)!}\\
& = & N_0t^{\mu-1}\sum^\infty_{r=0}(-1)^r\frac{(r+1)}{\Gamma(r\nu+\mu)}(ct)^{r\nu}\\
& = & N_0t^{\mu-1}\sum^\infty_{r=0}(-1)^r\frac{[\frac{1}{\nu}\left\{(r\nu+\mu-1)+(1+\nu-\mu)\right\}(ct)^{r\nu}}{\Gamma(r\nu+\mu)}\\
& = &N_0\frac{t^{\mu-1}}{\nu}\left[\sum^\infty_{r=0}\frac{(-c^\nu t^\nu)^r}{\Gamma(r\nu+\mu-1)}+(1+\nu-\mu\sum^\infty_{r=0}\frac{(-c^\nu t^\nu)^r}{\Gamma(r\nu+\mu)}\right]\\
& = & \mbox{R.H.S. of (30)}.
\end{eqnarray*}
\section {Conclusions}
The fractional kinetic equation (11) has been extended to generalized fractional equations (14), (21), (23), and (29). Their respective solutions are given in terms of the ordinary Mittag-Leffler function and their generalization, which can also be represented as FOX's H-functions. The ordinary and generalized Mittag-Leffler functions interpolate between a purely exponential law and power-like behavior of phenomena governed by ordinary kinetic equations and their fractional counterparts, respectively (Lang, 1999; Hilfer, 2000). A specific example for such behavior is the application of Tsallis statistics (Tsallis, 2002) to phenomena that may arise from fluctuations of temperature or energy dissipation rate (Lavagno and Quarati, 2002). The application of fractional kinetic equations to describe such phenomena has not been fully developed yet.   
\begin{center}
\bf{References}
\end{center}
\noindent
Agarwal, R.P.: 1953, A propos d'une note de M. Pierre Humbert, C.R. Acad.\par 
Sci. Paris {\bf 236}, 2031-2032.\\
Erd\'{e}lyi, A., Magnus, W., Oberhettinger, F., and Tricomi, F.G.: 1953,\par
Higher Transcendental Functions, Vol.{\bf 1}, McGraw-Hill,\par
New York-Toronto-London.\par
\smallskip
\noindent
Erd\'{e}lyi, A., Magnus, W., Oberhettinger, F., and Tricomi, F.G.: 1954, Tables\par 
of Integral Transforms, Vol.{\bf 1}, McGraw-Hill, New York-Toronto-London.\par
\smallskip
\noindent
Erd\'{e}lyi, A., Magnus, W., Oberhettinger, F., and Tricomi, F.G.: 1955, Higher\par 
Transcendental Functions, Vol.{\bf 3}, McGraw-Hill, New York-Toronto-London.\par
\smallskip
\noindent
Haubold, H.J. and Mathai, A.M.: 2000, The fractional kinetic equation and\par 
thermonuclear functions, Astrophysics and Space Science {\bf 327}, 53-63.\par
\smallskip
\noindent
Hilfer, R. (ed.): 2000, Applications of Fractional Calculus in Physics, World \par
Scientific, Singapore.\par
\smallskip
\noindent
Humbert, P.: 1953, Quelques resultats relatifs a'la fonction de Mittag-Leffler,\par 
C.R. Acad. Sci. Paris {\bf 236}, 1467-1468.\par
\smallskip
\noindent
Humbert, P. and Agarwal, R.P.: 1953, Sur la fonction de Mittag-Leffler et\par 
quelques-unes de ses generalisations, Bull. Sci. Math. (Ser.II) {\bf 77},\par 
180-185.\par
\smallskip
\noindent
Lang, K.R.: 1999, Astrophysical Formulae Vol. I (Radiation, Gas Processes\par
and High Energy Astrophysics) and Vol. II (Space, Time, Matter and\par
Cosmology), Springer-Verlag, Berlin-Heidelberg.\par
\smallskip
\noindent
Lavagno, A. and Quarati, P.: 2002, Classical and quantum non-extensive\par
statistics effects in nuclear many-body problems, Chaos, Solitons and\par
Fractals {\bf 13}, 569-580.\par
\smallskip
\noindent 
Miller, K.S. and Ross, B.: 1993, An Introduciton to the Fractional Calculus\par 
and Fractional Differential Equations, John Wiley and Sons, New York.\par
\smallskip
\noindent
Mittag-Leffler, G.M.: 1903, Sur la nouvelle fonction $E_\alpha(x)$, C.R. Acad. Sci.\par 
Paris, (Ser.II) {\bf 137}, 554-558.\par
\smallskip
\noindent
Mittag-Leffler, G.M.: 1905, Sur la representatin analytique d'une branche\par 
uniforme d'une fonciton monogene, Acta Math. {\bf 29}, 101-181.\par
\smallskip
\noindent
Oldham, K.B. and Spanier, J.: 1974, The Fractional Calculus: \par
Theory and Applications of Differentation and Integration to Arbitrary\par Order, 
Academic Press, New York.\par
\smallskip
\noindent
Srivastava, H.M. and Saxena, R.K.: 2001, Operators of fractional integration\par
and their applications, Applied Mathematics and Computation {\bf 118}, 1-52.\par
\smallskip
\noindent 
Tsallis, C.: 2002, Entropic nonextensivity: A possible measure of complexity,\par
Chaos, Solitons and Fractals {\bf 13}, 371-391.\par
\smallskip
\noindent 
Wiman, A.: 1905, Ueber den Fundamentalsatz in der Theorie der Funktionen\par
$E_\alpha (x)$, Acta Math. {\bf 29}, 191-201.\par
\smallskip
\noindent
Wiman, A.: 1905, Ueber die Nullstellen der Funktionen $E_\alpha(x)$, Acta Math.\par 
{\bf 29}, 217-234.
\end{document}